\documentclass{article}

\usepackage{amsmath,amssymb}

\newtheorem{theorem}{Theorem}[section]
\newtheorem{lemma}[theorem]{Lemma}
\newtheorem{remark}[theorem]{Remark}

\newtheorem{corollary}[theorem]{Corollary}
\newtheorem{example}[theorem]{Example}

\usepackage{color}

\newcommand{\qed}{$\hfill\square$\bigskip}
\newcommand{\g}[2]{\ensuremath{\langle #1,#2 \rangle}}

\newcommand{\proof}{\noindent{\bf Proof. }}
\newcommand{\ja}{\ensuremath{\mathcal{J}}}
\newcommand{\jw}{\ensuremath{\mathcal{J}_W}}

\title{Osserman and conformally Osserman manifolds with warped and twisted product structure}
\author{M. Brozos-V\' azquez$^*$, E. Garc\'ia-R\'io\thanks{Supported by projects MTM2006-01432 and
PGIDIT06PXIB207054PR.}, R. V\'azquez-Lorenzo$^*$}
\date{}

\begin{document}
\maketitle

\begin{abstract} We characterize Osserman and conformally Osserman Riemannian manifolds with the local structure of a warped product. By means of this approach we analyze the twisted product structure and obtain, as a consequence, that the only Osserman manifolds which
can be written as a twisted product are those of constant curvature.
\end{abstract}

\vspace{0.25cm}

\begin{enumerate}
\item[]{\underline{\normalsize\bf Mathematics Subject Classification.}}
53C20, 53C50
\item[]{\underline{\normalsize\bf Key words and phrases.}}
Jacobi operator, Osserman and conformally Osserman manifold, warped and twisted product structure.
\end{enumerate}

\section{Introduction}

The group of isometries of any Riemannian two-point-homogeneous
space acts transitively on the unit sphere bundle. This implies that
any Riemannian two-point-homogeneous manifold is Osserman, i.e.,
the eigenvalues of the Jacobi operator are constant on the unit
sphere bundle. Moreover, the converse was conjectured by Osserman
and proved to be true in dimension different from $16$
\cite{chi91,Nikp1,Nikp2}. Later, the Osserman concept was extended
to the conformal setting by analyzing the conformal Jacobi operator.
Remarkably both conditions are equivalent for Einstein manifolds of
dimension different from $4$, $16$. We refer to
\cite{blazic-gilkey-2004,blazic-gilkey-2005,blazic-gilkey-nikcevic-simon,
mer-prs} for basic results on conformally Osserman manifolds.

Warped products were introduced in \cite{bishop-oneill} as a tool to
construct Riemannian manifolds with non positive curvature.
Furthermore,  the fact that many cosmological models are warped
products reinforces their importance in Lorentzian signature.
Twisted products were introduced in $\cite{Ch}$ as a more general
structure, since they include a wider variety of metrics than warped
products. Both, warped and twisted products, are the underlying
structure of many geometrical situations and they play an important
role in conformal geometry.

In this paper we study the relation between these geometric
properties  and the algebraic structure of the metric tensor. More
specifically we study which Osserman manifolds are warped or twisted
products.  Our approach relies on the fact that the metric is
positive definite, and thus the results we obtain are in Riemannian
signature. As an exception, the analysis we carry out in
Section~\ref{Section-Dimension-four} does not require any
restriction on the signature, and thus we show that a necessary
condition for any four-dimensional warped product to be conformally
Osserman  is the local conformal flatness. In
Section~\ref{section-higher-dimension} we study Riemannian
conformally Osserman warped products in arbitrary dimension,
obtaining the following result.

\begin{theorem}\label{warped-conf-Osserman-conf-cha}
A Riemannian warped product $B\times_fF$ is conformally Osserman if and only if it is locally conformally flat.
\end{theorem}

Finally, we also give in Section~\ref{section-higher-dimension} a
complete characterization of Osserman twisted products for positive
definite metrics.

\begin{theorem}\label{theorem-Osserman-twisted}
A Riemannian twisted product $B\times_fF$ is pointwise Osserman if
and only if it is a space of constant sectional curvature.
\end{theorem}

\section{Preliminaries}

Let $(M,g)$ be a Riemannian manifold of dimension $n$. Let
$\nabla$ denote the Levi-Civita connection and $R$ the curvature
tensor chosen with the sign convention given by $R(X,Y)=\nabla_{[ X,Y ]}-[\nabla_X,\nabla_Y]$. For any orthonormal basis $\{E_1,\dots,E_n\}$, the Ricci tensor and the scalar curvature are defined as follows:
\[
\rho(X,Y)=\sum_{i=1}^n R(X,E_i,Y,E_i)\,,\qquad \tau=\sum_{i=1}^n \rho(E_i,E_i)\,.
\]
The curvature tensor $R$ decomposes as $R=C\bullet g+W$, where $C$ is the
Schouten tensor, $W$ is the Weyl tensor and $\bullet$ denotes de Kulkarni--Nomizu product.  Then let
\[
\ja(X)Y=R(X,Y)X, \quad \;\mbox{and} \;\quad \jw(X)Y=W(X,Y)X
\]
be the usual \emph{Jacobi operator} and the \emph{conformal Jacobi
operator}, respectively. $(M,g)$ is said to be \emph{pointwise
Osserman} (respectively, \emph{conformally Osserman}) if for every
$p\in M$ the Jacobi operator $\ja_p$ (respectively, the conformal
Jacobi operator) has constant eigenvalues on the unit tangent sphere
$S_p=\{x \in T_pM: g_p(x,x)=1\}$, where $T_pM$ is the space tangent
to $M$ on $p$. The relation between these two concepts was pointed
out in \cite{blazic-gilkey-nikcevic-simon}, showing that {\it
$(M,g)$ is pointwise Osserman if and only if $(M,g)$ is Einstein and
conformally Osserman}.

\medskip

\noindent{\bf  Warped product structures}

Let $(B,g_B)$ and $(F,g_F)$ be pseudo-Riemannian manifolds and let
$f:B\longrightarrow \mathbb{R}^+$ be a positive function on $B$. The
product manifold $M=B\times F$ endowed with the metric
\[
g=g_B\oplus f^2 g_F
\]
is called the \emph{warped product} $B\times_fF$ with \emph{base} $B$, \emph{fiber} $F$ and \emph{warping function} $f$.
The geometry of a warped product is described in terms of the curvature tensor as follows.

\begin{lemma}\label{lemma-curdef-warped}{\rm \cite{bishop-oneill}}
Let $(M,g)=B\times_f F$ be a warped product. Let $X, Y, Z\in
\mathfrak{X}(B)$ and let $U,V,W \in \mathfrak{X}(F)$. The
curvature tensor $R$ is given by:
\begin{enumerate}
\item[(i)] $R(X,Y)Z$ is the lift of $R^B(X,Y)Z$ on $B$,
\item[(ii)] $R(U,X)Y=\frac{H_f(X,Y)}{f}U$, \item[(iii)]
$R(X,Y)U=R(U,V)X=0$, \item[(iv)] $R(X,U)V=\frac{\g{U}{V}}{f}
\nabla_X(\nabla f)$, \item[(v)] $R(U,V)W=R^F(U,V)W-\frac{\g{\nabla
f}{\nabla f}}{f^2} \left(\g{U}{W}V-\g{V}{W}U \right)$.
\end{enumerate}
\end{lemma}

\medskip

\noindent{\bf Twisted product structures}

Let $(B,g_B)$ and $(F,g_F)$ be pseudo-Riemannian manifolds and let
$f:B\times F\longrightarrow \mathbb{R}^+$ be a positive function on
$B\times F$. The product manifold $M=B\times F$ endowed with the
metric
\[
g=g_B\oplus f^2 g_F
\]
is called the \emph{twisted product} $B\times_fF$. The same terminology of \emph{base} and \emph{fiber} applies in this case, whereas $f$ is referred to as the \emph{twisting function}. In the following lemma the geometry of a twisted
product is described. For simplicity in some of the expressions we use $\xi=Log\, f$
instead of $f$.

\begin{lemma}\label{lemma-twisted-cur}{\rm \cite{manolo-eduardo-kupeli-unal}, \cite{PR}}
Let $(M,g)=B\times_f F$ be a twisted product. For $X,Y,Z\in
\mathfrak{X}(B)$ and $U,V,W \in \mathfrak{X}(F)$ the curvature
tensor $R$ is given by:
\begin{enumerate}
\item[(i)]
$R(X,Y)Z= R^B(X,Y)Z$,
\item[(ii)]
$R(U,X)Y=
(H_\xi(X,Y)+X(\xi)Y(\xi))U$,
\item[(iii)]
$R(X,U)V=
g(U,V)(X(\xi)\nabla \xi+ h_\xi(X))- XV(\xi)U$, \item[(iv)]
$R(U,V)X= XV(\xi)U-XU(\xi)V$,
\item[(v)]
$\begin{array}[t]{@{}rcl}
                R(U,V)W=R^F(U,V)W&+& g(\nabla \xi,\nabla
                \xi)(g(U,W)V-g(V,W)U)\\
                \noalign{\smallskip}
                &+&H_\xi (V,W)U- H_\xi(U,W)V\\
                \noalign{\smallskip}
                &+&g(V,W)h_\xi(U)-g(U,W)h_\xi (V)\\
                \noalign{\smallskip}
                &+&W(\xi)(V(\xi)U-U(\xi)V)\\
                \noalign{\smallskip}
                &+&(U(\xi)g(V,W)- V(\xi)g(U,W))\nabla\xi\, .
           \end{array}$
\end{enumerate}
\end{lemma}

\begin{remark}\rm
The difference between warped and twisted products is pointed out by
reducing the base $B$ to a single point $\{p\}$. Then, warped
products $\{ p\}\times_fF$ become homothetic metrics on the fiber
$F$, while twisted products $\{ p\}\times_fF$ are in one to one
correspondence with conformal metrics on $F$.
\end{remark}

\section{Four-dimensional setting}\label{Section-Dimension-four}

Four-dimensional geometry is of special interest in the context of
Osserman and conformally Osserman geometries.  Three-dimensional
Osserman metrics are of constant curvature and moreover, the Weyl
tensor has no meaning in dimension three. Therefore, dimension four
is the first non-trivial case to be considered.

A proper feature of dimension four is that the Hodge  star operator
is idempotent and the curvature tensor further decomposes as
$R=C\bullet g+W_++W_-$, where $W_\pm$ denote the self-dual and
anti-self-dual components of the Weyl tensor. Moreover, since
\emph{a four-dimensional pseudo-Riemannian manifold is conformally
Osserman if and only if it is (anti-)self-dual} (see
\cite{blazic-gilkey-2005} and \cite{mer-prs} for a proof of this
fact), our approach in the present section relies on the study of
the self-dual and the anti-self-dual Weyl curvature operators.

\subsection{Four-dimensional warped products}

The following lemma, which relates self-duality and anti-self-duality, is a previous step to the characterization of conformally Osserman warped products in dimension four.

\begin{lemma}\label{lemma-warped-self=anti-selfdual}
A four-dimensional product manifold $B\times F$ is self-dual if and
only if it is anti-self-dual.
\end{lemma}
\proof Since any Lorentzian four-manifold is self-dual if and only
if it is anti-self-dual, next we consider  Riemannian $(++++)$ or
neutral signature $(++--)$ direct products $B\times F$. Let $p\in
B\times F$ be an arbitrary point. Let $\{e_1,e_2,e_3,e_4\}$ be an
orthonormal basis of the tangent space and let $\{e^1,e^2,e^3,e^4\}$
be the corresponding dual basis. Consider the following orthonormal
basis for the self-dual and anti-self-dual spaces:
\[
\begin{array}{c}
\Lambda^\pm={\rm Span\,}\{E_1^\pm=(e^1\wedge e^2\pm
\epsilon_3\epsilon_4 e^3\wedge e^4)/\sqrt{2},E_2^\pm=(e^1\wedge
e^3\mp
\epsilon_2\epsilon_4 e^2\wedge e^4)/\sqrt{2},\\
\noalign{\medskip} E_3^\pm=(e^1\wedge e^4\pm \epsilon_2\epsilon_3
e^2\wedge e^3)/\sqrt{2}\},
\end{array}
\]
where $\epsilon_i=g(e_i,e_i)$. The diagonal terms of the self-dual
and anti-self-dual matrix associated to these basis are given by:
\begin{equation}
\begin{array}{rcl}\label{diagonal-terms-weyl}
2W_{11}^+&=&W_{1212}+W_{3434}+2\epsilon_3\epsilon_4 W_{1234},\\
\noalign{\smallskip}
2W_{11}^-&=&W_{1212}+W_{3434}-2\epsilon_3\epsilon_4 W_{1234},\\
\noalign{\smallskip}
2W_{22}^+&=&W_{1313}+W_{2424}-2\epsilon_2\epsilon_4 W_{1324},\\
\noalign{\smallskip}
2W_{22}^-&=&W_{1313}+W_{2424}+2\epsilon_2\epsilon_4 W_{1324},\\
\noalign{\smallskip}
2W_{33}^+&=&W_{1414}+W_{2323}+2\epsilon_2\epsilon_3 W_{1423},\\
\noalign{\smallskip}
2W_{33}^-&=&W_{1414}+W_{2323}-2\epsilon_2\epsilon_3 W_{1423},
\end{array}
\end{equation}
but for a direct product we have $W_{1234}=W_{1324}=W_{1423}=0$, so
$W_{11}^+=W_{11}^-$, $W_{22}^+=W_{22}^-$ and $W_{33}^+=W_{33}^-$.
The remaining terms are:
\begin{equation}
\begin{array}{rcl}\label{non-diagonal-terms-weyl}
2W_{12}^+&=&W_{1213}-\epsilon_2\epsilon_4W_{1224}+\epsilon_3\epsilon_4W_{3413}-\epsilon_2\epsilon_3W_{3424},\\
\noalign{\smallskip}
2W_{12}^-&=&W_{1213}+\epsilon_2\epsilon_4W_{1224}-\epsilon_3\epsilon_4W_{3413}-\epsilon_2\epsilon_3W_{3424},\\
\noalign{\smallskip}
2W_{13}^+&=&W_{1214}+\epsilon_2\epsilon_3W_{1223}+\epsilon_3\epsilon_4W_{3414}+\epsilon_2\epsilon_4W_{3423},\\
\noalign{\smallskip}
2W_{13}^-&=&W_{1214}-\epsilon_2\epsilon_3W_{1223}-\epsilon_3\epsilon_4W_{3414}+\epsilon_2\epsilon_4W_{3423},\\
\noalign{\smallskip}
2W_{23}^+&=&W_{1314}+\epsilon_2\epsilon_3W_{1323}-\epsilon_2\epsilon_4W_{2414}-\epsilon_3\epsilon_4W_{2423},\\
\noalign{\smallskip}
2W_{23}^-&=&W_{1314}-\epsilon_2\epsilon_3W_{1323}+\epsilon_2\epsilon_4W_{2414}-\epsilon_3\epsilon_4W_{2423}.
\end{array}
\end{equation}

There are two different cases to be analyzed: $dim\,B=dim\,F=2$ and
$dim\,B=1$, $dim\,F=3$.

Suppose first that $dim\;B=2$, $dim\;F=2$. Let $e_1,e_2\in
\mathfrak{X}(B)$ and $e_3,e_4\in \mathfrak{X}(F)$. Then
\[
\begin{array}{c}
W_{1213}=W_{1224}=W_{3413}=W_{3424}=0,\\
\noalign{\smallskip}
W_{1214}=W_{1223}=W_{3414}=W_{3423}=0,\\
\noalign{\smallskip}
W_{1314}=-\frac{\epsilon_1}{2}\rho_{34},W_{1323}=-\frac{\epsilon_3}{2}\rho_{12},W_{2414}=-\frac{\epsilon_4}{2}\rho_{12},W_{2423}=-\frac{\epsilon_2}{2}\rho_{34}.
\end{array}
\]
Hence $W_{12}^+=W_{12}^-=W_{13}^+=W_{13}^-=0$. Now, since the
signature is Riemannian or neutral, we have
$\epsilon_1=\epsilon_2\epsilon_3\epsilon_4$, therefore
\[
W_{23}^+=-\frac14(\epsilon_1-\epsilon_2\epsilon_3\epsilon_4)\rho_{34}=0=W_{23}^-.
\]

Suppose now that $dim\;B=3$ and $dim\;F=1$. Let $e_1,e_2,e_3\in
\mathfrak{X}(B)$ and $e_4\in \mathfrak{X}(F)$. Then
\[
W_{1224}=W_{3413}=0, \quad W_{1214}=W_{3423}=0,\quad
W_{1314}=W_{2423}=0,
\]
so $W_{12}^+=W_{12}^-$, $W_{13}^+=-W_{13}^-$ and
$W_{23}^+=-W_{23}^-$.

The relations between the self-dual and the anti-self-dual
components show that in all cases the self-dual and anti-self-dual
operators simultaneously vanish. Hence the result
follows.\mbox{}\qed

The main result of this section reads as follows:

\begin{theorem}\label{theorem-conf-Osser-locconfflat-4dim}
A four-dimensional pseudo-Riemannian warped product manifold is conformally Osserman if and only if it is locally conformally flat.
\end{theorem}

\proof
If the signature is Lorentzian then the result follows from \cite{blazic-gilkey-nikcevic-simon}. Hence assume hereafter that the signature is Riemannian or neutral. A warped product $B\times_f F$ is in the conformal class of
$(B\times F,\frac{1}{f^2}g_B\oplus g_F)$, which is a direct product.
Thus, since the two properties under consideration are conformal
invariants, it suffices for our purpose to restrict our analysis
to direct products. From Lemma~\ref{lemma-warped-self=anti-selfdual} a direct product is self-dual if and only if it is anti-self-dual and the result follows.
\qed

\subsection{Four-dimensional twisted products}

Note that the relations given in the proof of
Lemma~\ref{lemma-warped-self=anti-selfdual} between the self-dual
and anti-self-dual components still hold if we consider a twisted
product $B\times_f F$ with $dim\; F=1$. Hence
Lemma~\ref{lemma-warped-self=anti-selfdual} also holds in this
particular case. Moreover, for a twisted product $B\times_f F$ with
$dim\; B=1$ we consider a conformal change by $\frac{1}{f^2}$ to get
$(B\times F,\frac{1}{f^2}g_B\oplus g_F)$, so that the dimension of
the fiber is $1$. Again the result holds since the conditions under
consideration are conformally invariant, thus we obtain the
following result.

\begin{theorem}\label{theorem-4dim-twisted}
Let $B\times_f F$ be a four-dimensional pseudo-Riemannian twisted
product with $dim\,B=1$ or $dim\,B=3$. Then $(M,g)$ is conformally
Osserman if and only if it is locally conformally flat.
\end{theorem}

In view of Theorem~\ref{theorem-4dim-twisted} one may wonder if this
result holds in general for a four-dimensional twisted product. This
is not the case, as next example shows.
\begin{example}\rm\label{example-twisted-nowarped}
Consider the twisted product $ (M,g)=\mathbb{R}^2\times_f
\mathbb{R}^2$, \,with twisting function\, $\displaystyle
f(x_1,x_2,x_3,x_4)=e^{x_1 x_3-x_2 x_4}$. A lengthy but
straightforward calculation shows that
\[
W^+=\left(\begin{array}{ccc}
0&0&\frac{1}2(1+e^{x_1 x_3-x_2
x_4})\\
0&0&0\\
\frac{1}2(1+e^{x_1 x_3-x_2 x_4})&0&0
\end{array}\right) \qquad\mbox{and}\qquad W^-=0.
\]
Hence $(M,g)$ is self-dual but it is not anti-self-dual and,
consequently, it is conformally Osserman but not locally
conformally flat.
\end{example}

Although we have just given a general negative answer to the extension of Theorem~\ref{theorem-conf-Osser-locconfflat-4dim} to twisted products, the following result shows the non-existence of compact
(anti-)self-dual twisted products which are not locally
conformally flat if the metric is positive definite.

\begin{theorem}
Let $(M,g)=B\times_f F$ be a compact Riemannian twisted product such
that $dim\;B,F=2$. Then $B\times_f F$ is conformally Osserman if and
only if it is locally conformally flat. Moreover, in such a case it
is actually a warped product.
\end{theorem}

\proof Since $B\times_f F$ is conformally Osserman, it is self-dual
or anti-self-dual. Suppose $B\times_fF$ is self-dual (reversing the
orientation interchanges the roles of the self-dual and  the
anti-self-dual components). Since $B$ and $F$ are $2$-dimensional
and oriented, let $J^B$ and $J^F$ be orthogonal complex structures
on $B$ and $F$, respectively. Then  $B\times_fF$ is Hermitian and
opposite Hermitian just considering the complex structure
$J^B\oplus\pm J^F$. If  $B\times_fF$ is self-dual, then results in
\cite{apostolov-davidov-muskarov} show that  $B\times_fF$ is locally
conformally flat or conformally equivalent to $\mathbb{CP}^2$ with
the Fubini-Study metric or to a compact quotient of the unit ball in
$\mathbb{C}^2$ with the Bergman metric. Hence, it follows that
$B\times_fF$ is locally conformally flat or a locally conformally
K\"{a}hler manifold. Assume $J$ is locally conformally K\"{a}hler. Let $X$
be a vector on the base and $U$ a vector on the fiber. Hence
$\{X,JX,U,JU\}$ is an adapted basis for the product. Since
$(B\times_f F,J)$ is locally conformally K\"{a}hler, then there exists
$\psi$ such that $(B\times F, \psi\, g^B\oplus \psi\, f^2 g^F)$ is
K\"{a}hler on a suitable open set (note that this is a doubly twisted
product). Then
\[
\begin{array}{rcl}
(\nabla_X J)U&=&\nabla_X (JU)-J\nabla_X U\\
\noalign{\smallskip} &=&JU(\ln\, \psi)X+X(\ln\, \psi\,
f^2)JU-U(\ln\,
\psi)JX-X(\ln\, \psi\, f^2)JU\\
\noalign{\smallskip}
&=&JU(\ln\, \psi)X-U(\ln\, \psi)JX,\\
\noalign{\bigskip}
(\nabla_U J)X&=&\nabla_U (JX)-J\nabla_U X\\
\noalign{\smallskip}
&=&JX(\ln\, \psi\, f^2)U+U(\ln\, \psi)JX-X(\ln\, \psi\, f^2)JU-U(\ln\, \psi)JX\\
\noalign{\smallskip} &=&JX(\ln\, \psi\, f^2)U-X(\ln\, \psi\, f^2)JU,
\end{array}
\]
from where $U(\ln\, \psi)=0$, $JU(\ln\, \psi)=0$, $X(\ln\, \psi
f^2)=0$ and $JX(\ln\, \psi f^2)=0$. This implies that $\psi$ is
constant over $F$ and one proceeds in an analogous way to show that
$\psi f^2$ is constant over $B$ too. Therefore $f$ decomposes as a
product $f=f_B f_F$ so that $f_B$ is constant on $F$ and $f_F$ is
constant on $B$, hence $B\times_f F$ is indeed a warped product and
by Theorem~\ref{theorem-conf-Osser-locconfflat-4dim} it is locally
conformally flat. \qed

\section{Higher dimensional setting}\label{section-higher-dimension}

Motivated by the results of previous section, one may wonder if an
arbitrary conformally Osserman warped product is locally conformally
flat. It has been shown in \cite{blazic-gilkey-nikcevic-simon} that
any conformally  Osserman Lorentzian manifold is locally conformally
flat but the following is a counterexample to this fact in arbitrary
indefinite signature (not Lorentzian).

\begin{example}\rm
Consider the manifold $\mathcal{B}\times \mathbb{R}^k_\nu$, where
$\mathcal{B}=(\mathbb{R}^4,g_\mathcal{B})$ and $g_B$ is given by
\[
g_\mathcal{B}(x_1,x_2,x_3,x_4)=\left(\begin{array}{cccc}
                    0&0&1&0\\
                    0&0&0&1\\
                    1&0&0&c(x_3,x_4)\\
                    0&1&c(x_3,x_4)&0
                    \end{array}\right),
\]
where $(x_1,x_2,x_3,x_4)$ are coordinates in $\mathcal{B}$. Note
that $\mathcal{B}=(\mathbb{R}^4,g_\mathcal{B})$ is a strictly Walker
metric with signature $(2,2)$ and $\mathbb{R}_\nu^k$ is endowed with
an Euclidian metric of signature $(\nu,k-\nu)$. Then the only
non-null component of the curvature, up to the usual symmetries, is
\[
R(\partial_3,\partial_4,\partial_3,\partial_4)=\partial_3\partial_4c(x_3,x_4).
\]
A straightforward calculation shows that $M$ is Ricci flat and thus the only non-null component of the Weyl tensor, up to the usual symmetries, is
\[
W(\partial_3,\partial_4,\partial_3,\partial_4)=\partial_3\partial_4c(x_3,x_4).
\]
Then $\mathcal{B}\times \mathbb{R}^k_\nu$ is conformally Osserman
with nilpotent conformal Jacobi operator; however it is not locally
conformally flat unless $\partial_3\partial_4c(x_3,x_4)=0$.
\end{example}

In what follows we show that the result in theorems
\ref{theorem-conf-Osser-locconfflat-4dim} and
\ref{theorem-4dim-twisted} still holds in the Riemannian setting for
arbitrary dimension. We proceed in an analogous way to the previous
section and, as a preliminary step, we characterize conformally
Osserman direct products.

\begin{lemma}\label{clasificacion-produto-conf-Osserman}
Let $(M,g)$ be a conformally Osserman Riemannian
manifold which decomposes as a direct product $B\times F$. Then $W^M=0$.
\end{lemma}
\proof Let $p\in M$ be an arbitrary point. To establish notation,
set $b=dim\,B$ and $d=dim\,F$ and use superindexes $B$ and $F$ to
refer to manifolds $B$ and $F$, respectively. On $T_pM$ we can
choose an orthonormal basis $\{e_1,\dots,e_b,f_1,\dots,f_d\}$, with
$\{e_i\}\subset T_p^BM$ and $\{f_i\}\subset T_p^FM$, which
diagonalizes the Ricci tensor.

Note that $\mathcal{J}_W(e_i)e_j,\mathcal{J}_W(f_i)e_j\in T_p^BM$
and $\mathcal{J}_W(f_i)f_j,\mathcal{J}_W(e_i)f_j\in T_p^FM$.
Indeed, from the expressions of the curvature in
Lemma~\ref{lemma-curdef-warped} we compute:
\begin{equation}\label{formula-autovalores}
\begin{array}{rcl}
\mathcal{J}_W(e_i)e_j&=&\mathcal{J}_R(e_i)e_j-\frac{1}{n-2}\left(\rho(e_i,e_i)+\rho(e_j,e_j)-\frac{\tau}{n-1}\right)e_j,\,\mbox{and} \\
\noalign{\medskip}
\mathcal{J}_W(e_i)f_j&=&-\frac{1}{n-2}\left(\rho(e_i,e_i)+\rho(f_j,f_j)-\frac{\tau}{n-1}\right)f_j.
\end{array}
\end{equation}
Therefore, $f_j$ is an eigenvector for every
$\mathcal{J}_W(e_i)$ (analogously, $e_i$ is an eigenvector for
every $\mathcal{J}_W(f_j)$). Also notice that mixed terms of the
Weyl tensor vanish, that is, $W(e_i,f_j)e_k=0$ and
$W(f_i,e_j)f_k=0$ if $i\neq k$.

Suppose $\mathcal{J}_W(e_1)f_j=\lambda\,f_j$,
$\mathcal{J}_W(e_1)f_k=\mu\, f_k$ with $j\neq k$. Recall the Raki\'c
duality principle \cite{rakic-duality}: \emph{``for an Osserman
Riemannian manifold $\ja_A(X)Y=\lambda\, Y$ if and only if
$\ja_A(Y)X=\lambda\, X$"}. This principle also applies for
conformally Osserman Riemannian manifolds (indeed it is true in a
purely algebraic context for any Osserman algebraic curvature
tensor, see \cite{G02}). Hence we compute
\[
\begin{array}{rcl}
\mathcal{J}_W(\cos\theta f_j+\sin\theta f_k)e_1&=&\cos^2\theta \mathcal{J}_W(f_j)e_1+\sin^2\theta \mathcal{J}_W(f_k)e_1\\
\noalign{\smallskip}
&&+\cos\theta \sin\theta (W(f_j,e_1)f_k+W(f_k,e_1)f_j)\\
\noalign{\medskip}
&=&\cos^2\theta \mathcal{J}_W(f_j)e_1+\sin^2\theta \mathcal{J}_W(f_k)e_1\\
\noalign{\medskip}
&=&\left(\cos^2\theta\lambda+\sin^2\theta\mu\right)e_1.
\end{array}
\]
This shows that $e_1$ is an eigenvector for $\mathcal{J}_W(\cos\theta
f_j+\sin\theta f_k)$ associated to the eigenvalue
$\cos^2\theta\lambda+\sin^2\theta\mu$; but, since the
eigenvalues are constant, we conclude that $\lambda=\mu$. By
repeating this argument, we show that all the eigenvalues of
$\mathcal{J}_W(e_1)$ associated to eigenvectors $f_j$ in $T_p^FM$
are equal.

Next, we show that all the eigenvalues of
$\mathcal{J}_W(e_1)$ are indeed equal. Take a unitary vector $x\in T_p^BM$ such that
$\mathcal{J}_W(e_1)x=\nu x$. Then
\[
\begin{array}{rcl}
\mathcal{J}_W(\cos\theta e_1+\sin\theta f_1)x&=&\cos^2\theta \mathcal{J}_W(e_1)x+\sin^2\theta \mathcal{J}_W(f_1)x\\
\noalign{\smallskip}
&&+\cos\theta\sin\theta\left(W(e_1,x)f_1+W(f_1,x)e_1\right)\\
\noalign{\medskip}
&=&\left(\cos^2\theta\nu+\sin^2\theta\lambda\right)x,
\end{array}
\]
and necessarily $\lambda=\nu$. Since the trace of $\mathcal{J}_W(\cdot)$ is zero, we conclude
that all the eigenvalues of $\mathcal{J}_W(\cdot)$ vanish and hence
the Weyl tensor is zero. \qed

\noindent{\bf Proof of Theorem~\ref{warped-conf-Osserman-conf-cha}.}
It follows immediately from Lemma \ref{clasificacion-produto-conf-Osserman} just using that any warped product metric is in the conformal class of a product metric.\qed

\subsection{Osserman condition on twisted products}

It was shown in \cite{miguel-eduardo-ramon} that a locally
conformally flat twisted product with factors of dimension greater
than one may be expressed as a warped product. Hence it arises as a
natural question whether a conformally Osserman twisted product can
be reduced to a warped one. The answer to this question is positive
for odd dimensions, since in this case the conformally Osserman
condition is equivalent to local conformal flatness (see
\cite{blazic-gilkey-nikcevic-simon}). On the contrary conformally
Osserman twisted products with base and fiber of dimension $\geq 2$
do not necessarily reduce to warped products as
Example~\ref{example-twisted-nowarped} shows. However, we see in
this section that manifolds with twisted product structure and fiber
or base of dimension one behave quite differently.


\begin{theorem}\label{twisted-dim1-conf-cha}
A twisted product $B\times_f F$, with $dim\;B=1$ or $dim\,F=1$, is
conformally Osserman if and only if it is locally conformally flat.
\end{theorem}

\proof If $dim\;B=1$, $B\times_f F$ is in the conformal class of
$F\times_{1/f} B$ whose fiber is one-dimensional, so we may suppose
without loss of generality that $dim\; F=1$. Let $p\in M$ be an
arbitrary point. Let $v$ be a unit vector on the fiber $F$. Now
consider $\mathcal{J}_W(v)$ and take a basis of orthonormal
eigenvectors $\{e_1,\dots,e_{n-1}\}$ tangent to the base of the
twisted product. Denote by $\{\lambda_1,\dots,\lambda_{n-1},0\}$ the
eigenvalues associated to the orthonormal basis
$\{e_1,\dots,e_{n-1},v\}$. Note that for $i\neq j$,
\[
W(e_i,v,e_j,v)=\langle \mathcal{J}_W(v)e_i,e_j\rangle=0,
\]
and, for any $i,j,k$ with $i\neq j$, $j\neq k$, $k\neq i$,
\[
\begin{array}{rcl}
W(e_i,v,e_j,e_k)&=&W(e_j,e_k,e_i,v)\\
\noalign{\medskip}
&=&R(e_j,e_k,e_i,v)-\frac{1}{n-2}\left(\rho(e_j,e_i)\langle
e_k,v\rangle+\rho(e_k,v)\langle e_j,e_i\rangle\right.\\
\noalign{\smallskip} &&\qquad\left.-\rho(e_j,v)\langle
e_k,e_i\rangle-\rho(e_k,e_i)\langle
e_j,v\rangle\right.\\
\noalign{\smallskip} &&\qquad\left.-\frac{\tau}{n-1}\langle
e_j,e_i\rangle \langle e_k,v\rangle-\langle e_j,v\rangle \langle
e_i,e_k\rangle\right)\\
&=&0.
\end{array}
\]
Also note that $W(e_i,v,e_j,e_i)=0$. Hence $W(e_i,v)e_j=0$ for
$i\neq j$ and $\mathcal{J}_W(e_i)v=\lambda_i v$. Then
\[
\begin{array}{rcl}
\mathcal{J}_W(\cos\theta e_i+\sin\theta e_j)v&=&\cos^2\theta
\mathcal{J}_W(e_i)v+\sin^2\theta \mathcal{J}_W(e_j)v\\
\noalign{\smallskip}
&&+\cos\theta\sin\theta(W(e_i,v)e_j+W(e_j,v)e_i)\\
\noalign{\medskip} &=&(\cos^2\theta \lambda_i +\sin^2\theta
\lambda_j) v.
\end{array}
\]
Since $M$ is conformally Osserman, the eigenvalues of the
conformal Jacobi operator are constant. Therefore
$\lambda_i=\lambda_j$. But, since ${\rm
tr}(\mathcal{J}_W(\cdot))=0$, all the eigenvalues are zero and
thus the Weyl tensor vanishes. Hence $M$ is locally conformally
flat. \qed

Finally we take advantage of
Theorem~\ref{twisted-dim1-conf-cha} to show that every Osserman twisted product has constant sectional curvature in the Riemannian setting.
\medskip

\noindent{\bf Proof of Theorem~\ref{theorem-Osserman-twisted}.}
Assume $B\times_fF$ is Osserman. Then it is Einstein and conformally
Osserman \cite{blazic-gilkey-nikcevic-simon}. On the one hand, if
$dim\,F>1$ then, by results in \cite{manolo-eduardo-kupeli-unal},
$B\times_fF$ is indeed a warped product and, from
Theorem~\ref{warped-conf-Osserman-conf-cha}, $B\times_fF$ is locally
conformally flat. On the other hand, if $dim\,F=1$ we conclude from
Theorem~\ref{twisted-dim1-conf-cha} that $B\times_fF$ is locally
conformally flat too. Finally, since $B\times_fF$ is Einstein and
locally conformally flat, it has constant sectional curvature. \qed

Since any two-point-homogeneous manifold is Osserman, the following
result is an immediate consequence of
Theorem~\ref{theorem-Osserman-twisted}.

\begin{corollary}\label{cor:osserman-twisted}
A two-point homogeneous space can be decomposed as a twisted product if and only if it is of constant sectional curvature.
\end{corollary}

\bigskip
\noindent Authors' addresses:
\medskip

\noindent
Miguel Brozos-V\'{a}zquez, Eduardo Garc\'{\i}a-R\'{\i}o, Ram\'{o}n V\'{a}zquez-Lorenzo

\noindent Faculty of Mathematics, University of Santiago de
Compostela, 15782 Santiago de Compostela, Spain,

\noindent (mbrozos@edu.xunta.es, xtedugr@usc.es, ravazlor@edu.xunta.es)

\end{document}